# DISCUSSION OF "LEAST ANGLE REGRESSION" BY EFRON ET AL.

By Berwin A. Turlach

*University of Western Australia*

I would like to begin by congratulating the authors (referred to below as EHJT) for their interesting paper in which they propose a new variable selection method (LARS) for building linear models and show how their new method relates to other methods that have been proposed recently. I found the paper to be very stimulating and found the additional insight that it provides about the Lasso technique to be of particular interest.

My comments center around the question of how we can select linear models that conform with the marginality principle [Nelder (1977, 1994) and McCullagh and Nelder (1989)]; that is, the response surface is invariant under scaling and translation of the explanatory variables in the model. Recently one of my interests was to explore whether the Lasso technique or the nonnegative garrote [Breiman (1995)] could be modified such that it incorporates the marginality principle. However, it does not seem to be a trivial matter to change the criteria that these techniques minimize in such a way that the marginality principle is incorporated in a satisfactory manner.

On the other hand, it seems to be straightforward to modify the LARS technique to incorporate this principle. In their paper, EHJT address this issue somewhat in passing when they suggest toward the end of Section 3 that one first fit main effects only and interactions in a second step to control the order in which variables are allowed to enter the model. However, such a two-step procedure may have a somewhat less than optimal behavior as the following, admittedly artificial, example shows.

Assume we have a vector of explanatory variables $X = (X_1, X_2, \ldots, X_{10})$ where the components are independent of each other and $X_i$, $i = 1, \ldots, 10$, follows a uniform distribution on $[0, 1]$. Take as model

$$Y = (X_1 - 0.5)^2 + X_2 + X_3 + X_4 + X_5 + \varepsilon, \tag{1}$$

where $\varepsilon$ has mean zero, has standard deviation 0.05 and is independent of $X$.

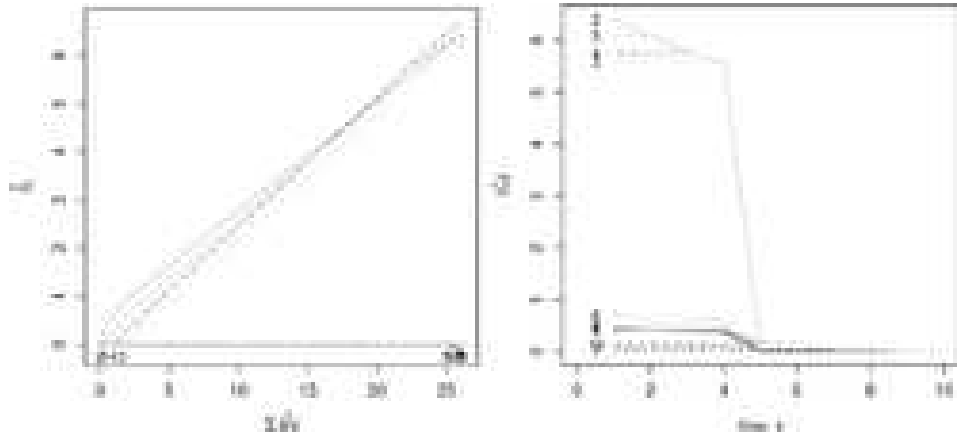

Fig. 1. *LARS analysis of simulated data with main terms only:* (left) *estimates of regression coefficients $\hat{\beta}_j$, $j = 1, \ldots, 10$, plotted versus $\sum |\hat{\beta}_j|$;* (right) *absolute current correlations as functions of LARS step.*

It is not difficult to verify that in this model $X_1$ and $Y$ are uncorrelated. Moreover, since the $X_i$'s are independent, $X_1$ is also uncorrelated with any residual vector coming from a linear model formed only by explanatory variables selected from $\{X_2, \ldots, X_{10}\}$.

Thus, if one fits a main effects only model, one would expect that the LARS algorithm has problems identifying that $X_1$ should be part of the model. That this is indeed the case is shown in Figure 1. The picture presents the result of the LARS analysis for a typical data set generated from model (1); the sample size was $n = 500$. Note that, unlike Figure 3 in EHJT, Figure 1 (and similar figures below) has been produced using the standardized explanatory variables and no back-transformation to the original scale was done.

For this realization, the variables are selected in the sequence $X_2$, $X_5$, $X_4$, $X_3$, $X_6$, $X_{10}$, $X_7$, $X_8$, $X_9$ and, finally, $X_1$. Thus, as expected, the LARS algorithm has problems identifying $X_1$ as part of the model. To further verify this, 1000 different data sets, each of size $n = 500$, were simulated from model (1) and a LARS analysis performed on each of them. For each of the 10 explanatory variables the step at which it was selected to enter the model was recorded. Figure 2 shows for each of the variables the (percentage) histogram of these data.

It is clear that the LARS algorithm has no problems identifying that $X_2, \ldots, X_5$ should be in the model. These variables are all selected in the first four steps and, not surprisingly given the model, with more or less equal probability at any of these steps. $X_1$ has a chance of approximately 25% of being selected as the fifth variable, otherwise it may enter the model at step 6, 7, 8, 9 or 10 (each with probability roughly 15%). Finally, each of the variables $X_6$ to $X_{10}$ seems to be selected with equal probability anywhere between step 5 and step 10.

This example shows that a main effects first LARS analysis followed by a check for interaction terms would not work in such cases. In most cases the LARS analysis would miss $X_1$ as fifth variable and even in the cases where it was selected at step 5 it would probably be deemed to be unimportant and excluded from further analysis.

How does LARS perform if one uses from the beginning all 10 main effects and all 55 interaction terms? Figure 3 shows the LARS analysis for the same data used to produce Figure 1 but this time the design matrix was augmented to contain all main effects and all interactions. The order in which the variables enter the model is $X_{2:5} = X_2 \times X_5$, $X_{2:4}$, $X_{3:4}$, $X_{2:3}$, $X_{3:5}$, $X_{4:5}$, $X_{5:5} = X_5^2$, $X_4$, $X_3$, $X_2$, $X_5$, $X_{4:4}$, $X_{1:1}$, $X_{1:6}$, $X_{1:9}$, $X_1$, .... In this example the last of the six terms that are actually in model (1) was selected by the LARS algorithm in step 16.

Using the same 1000 samples of size $n = 500$ as above and performing a LARS analysis on them using a design matrix with all main and interaction terms shows a surprising result. Again, for each replication the step at which a variable was selected into the model by LARS was recorded and Figure 4 shows for each variable histograms for these data. To avoid cluttering, the histograms in Figure 4 were truncated to $[1, 20]$; the complete histograms are shown on the left in Figure 7.

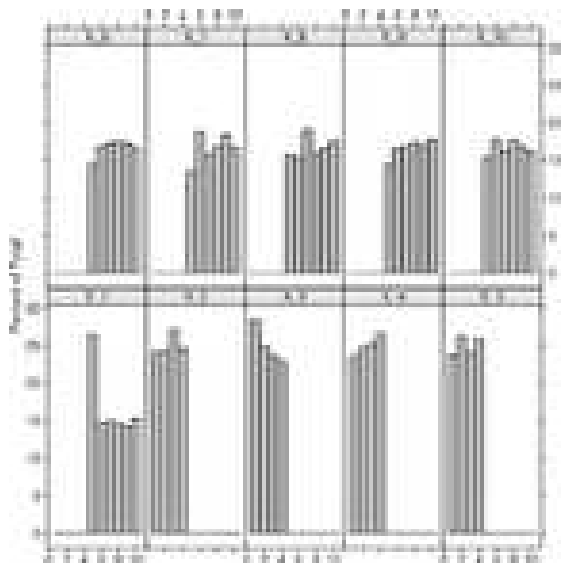

FIG. 2. *Percentage histogram of step at which each variable is selected based on* 1000 *replications: results shown for LARS analysis using main terms only.*

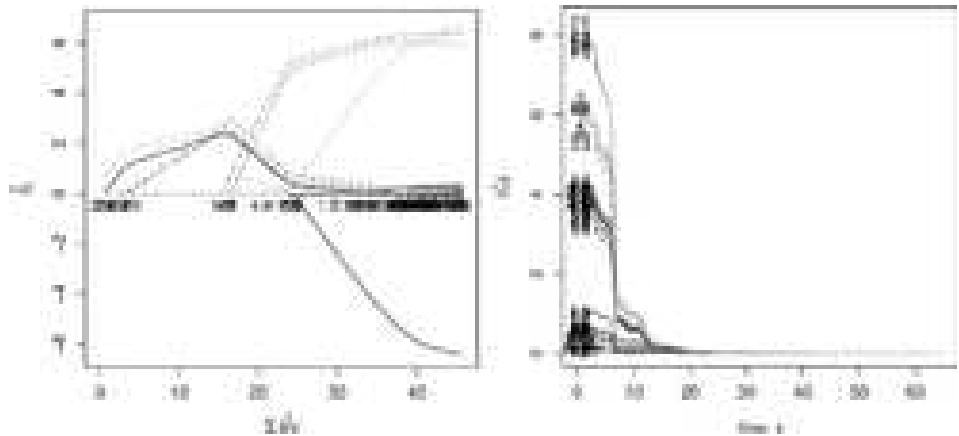

Fig. 3. *LARS analysis of simulated data with main terms and interaction terms:* (left) *estimates of regression coefficients* $\hat{\beta}_j$, $j = 1, \ldots, 65$, *plotted versus* $\sum |\hat{\beta}_j|$; (right) *absolute current correlations as functions of LARS step.*

The most striking feature of these histograms is that the six interaction terms $X_{i:j}$, $i, j \in \{2, 3, 4, 5\}$, $i < j$, were always selected first. In no replication was any of these terms selected after step 6 and no other variable was ever selected in the first six steps. That one of these terms is selected as the first term is not surprising as these variables have the highest correlation with the response variable $Y$. It can be shown that the covariance of these interaction terms with $Y$ is by a factor $\sqrt{12/7} \approx 1.3$ larger than the covariance between $X_i$ and $Y$ for $i = 2, \ldots, 5$. But that these six interaction terms dominate the early variable selection steps in such a manner came as a bit as a surprise.

After selecting these six interaction terms, the LARS algorithm then seems to select mostly $X_2$, $X_3$, $X_4$ and $X_5$, followed soon by $X_{1:1}$ and $X_1$. However, especially the latter one seems to be selected rather late and other terms may be selected earlier. Other remarkable features in Figure 4 are the peaks in histograms of $X_{i:i}$ for $i = 2, 3, 4, 5$; each of these terms seems to have a fair chance of being selected before the corresponding main term and before $X_{1:1}$ and $X_1$.

One of the problems seems to be the large number of interaction terms that the LARS algorithm selects without putting the corresponding main terms into the model too. This behavior violates the marginality principle. Also, for this model, one would expect that ensuring that for each higher-order term the corresponding lower-order terms are in the model too would alleviate the problem that the six interaction terms $X_{i:j}$, $i, j \in \{2, 3, 4, 5\}$, $i < j$, are always selected first.

I give an alternative description of the LARS algorithm first before I show how it can be modified to incorporate the marginality principle. This description is based on the discussion in EHJT and shown in Algorithm 1.

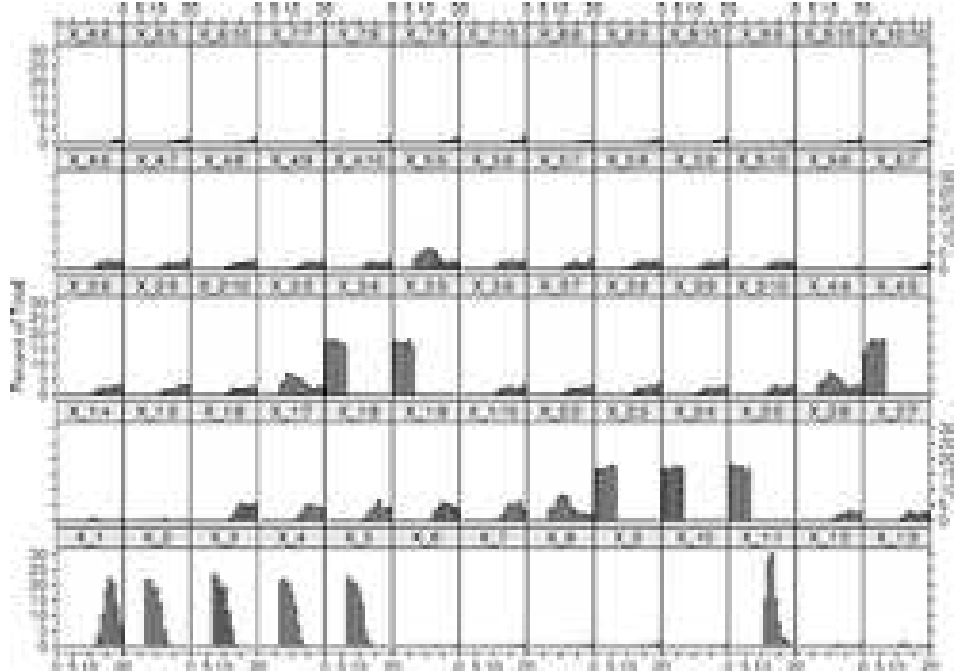

FIG. 4. *Percentage histogram of step at which each variable is selected based on* 1000 *replications: results shown for variables selected in the first* 20 *steps of a LARS analysis using main and interaction terms.*

ALGORITHM 1 (An alternative description of the LARS algorithm).

1. Set $\hat{\boldsymbol{\mu}}_0 = \mathbf{0}$ and $k = 0$.
2. **repeat**
3.    Calculate $\hat{\mathbf{c}} = X'(\mathbf{y} - \hat{\boldsymbol{\mu}}_k)$ and set $\hat{C} = \max_j\{|\hat{c}_j|\}$.
4.    Let $\mathcal{A} = \{j : |\hat{c}_j| = \hat{C}\}$.
5.    Set $X_\mathcal{A} = (\cdots \mathbf{x}_j \cdots)_{j \in \mathcal{A}}$ for calculating $\bar{\mathbf{y}}_{k+1} = (X'_\mathcal{A} X_\mathcal{A})^{-1} X'_\mathcal{A} \mathbf{y}$ and $\mathbf{a} = X'_\mathcal{A}(\bar{\mathbf{y}}_{k+1} - \hat{\boldsymbol{\mu}}_k)$.
6.    Set

$$\hat{\boldsymbol{\mu}}_{k+1} = \hat{\boldsymbol{\mu}}_k + \hat{\gamma}(\bar{\mathbf{y}}_{k+1} - \hat{\boldsymbol{\mu}}_k),$$

   where, if $\mathcal{A}^c \neq \varnothing$,

$$\hat{\gamma} = \min_{j \in \mathcal{A}^c}{}^+ \left\{ \frac{\hat{C} - \hat{c}_j}{\hat{C} - a_j}, \frac{\hat{C} + \hat{c}_j}{\hat{C} + a_j} \right\},$$

   otherwise set $\hat{\gamma} = 1$.
7.    $k \leftarrow k + 1$.
8. **until** $\mathcal{A}^c = \varnothing$.

We start with an estimated response $\hat{\boldsymbol{\mu}}_0 = \mathbf{0}$ and then iterate until all variables have been selected. In each iteration, we first determine (up to a constant factor) the correlation between all variables and the current residual vector. All variables whose absolute correlation with the residual vector equals the maximal achievable absolute correlation are chosen to be in the model and we calculate the least squares regression response, say $\bar{\mathbf{y}}_{k+1}$, using these variables. Then we move from our current estimated response $\hat{\boldsymbol{\mu}}_k$ toward $\bar{\mathbf{y}}_{k+1}$ until a new variable would enter the model.

Using this description of the LARS algorithm, it seems obvious how to modify the algorithm such that it respects the marginality principle. Assume that for each column $i$ of the design matrix we set $d_{ij} = 1$ if column $j$ should be in the model whenever column $i$ is in the model and zero otherwise; here $j \neq i$ takes values in $\{1, \ldots, m\}$, where $m$ is the number of columns of the design matrix. For example, abusing this notation slightly, for model (1) we might set $d_{1:1,1} = 1$ and all other $d_{1:1,j} = 0$; or $d_{1:2,1} = 1$, $d_{1:2,2} = 1$ and all other $d_{1:2,j}$ equal to zero.

Having defined such a dependency structure between the columns of the design matrix, the obvious modification of the LARS algorithm is that when adding, say, column $i$ to the selected columns one also adds all those columns for which $d_{ij} = 1$. This modification is described in Algorithm 2.

ALGORITHM 2 (The modified LARS algorithm).

1. Set $\hat{\boldsymbol{\mu}}_0 = \mathbf{0}$ and $k = 0$.
2. **repeat**
3.     Calculate $\hat{\mathbf{c}} = X'(\mathbf{y} - \hat{\boldsymbol{\mu}}_k)$ and set $\hat{C} = \max_j \{|\hat{c}_j|\}$.
4.     Let $\mathcal{A}_0 = \{j : |\hat{c}_j| = \hat{C}\}$, $\mathcal{A}_1 = \{j : d_{ij} \neq 0, i \in \mathcal{A}_0\}$ and $\mathcal{A} = \mathcal{A}_0 \cup \mathcal{A}_1$.
5.     Set $X_\mathcal{A} = (\cdots \mathbf{x}_j \cdots)_{j \in \mathcal{A}}$ for calculating $\bar{\mathbf{y}}_{k+1} = (X'_\mathcal{A} X_\mathcal{A})^{-1} X'_\mathcal{A} \mathbf{y}$ and $\mathbf{a} = X'_\mathcal{A}(\bar{\mathbf{y}}_{k+1} - \hat{\boldsymbol{\mu}}_k)$.
6.     Set
$$\hat{\boldsymbol{\mu}}_{k+1} = \hat{\boldsymbol{\mu}}_k + \hat{\gamma}(\bar{\mathbf{y}}_{k+1} - \hat{\boldsymbol{\mu}}_k),$$
where, if $\mathcal{A}^c \neq \varnothing$,
$$\hat{\gamma} = \min_{j \in \mathcal{A}^c}{}^+ \left\{ \frac{\hat{C} - \hat{c}_j}{\hat{C} - a_j}, \frac{\hat{C} + \hat{c}_j}{\hat{C} + a_j} \right\},$$
otherwise set $\hat{\gamma} = 1$.
7.     $k \leftarrow k + 1$.
8. **until** $\mathcal{A}^c = \varnothing$.

Note that compared with the original Algorithm 1 only the fourth line changes. Furthermore, for all $i \in \mathcal{A}$ it is obvious that for $0 \leq \gamma \leq 1$ we have

(2) $$|\hat{c}_i(\gamma)| = (1 - \gamma)|\hat{c}_i|,$$

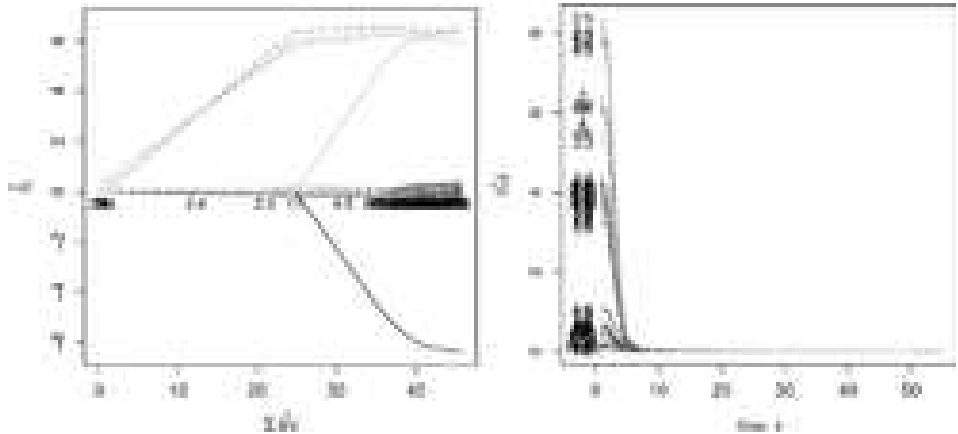

FIG. 5. *Modified LARS analysis of simulated data with main terms and interaction terms:* (left) *estimates of regression coefficients* $\hat{\beta}_j$, $j = 1, \ldots, 65$, *plotted versus* $\sum |\hat{\beta}_j|$; (right) *absolute current correlations as functions of* $k = \#\mathcal{A}^c$.

where $\hat{\mathbf{c}}(\gamma) = X'(\mathbf{y} - \hat{\boldsymbol{\mu}}(\gamma))$ and $\hat{\boldsymbol{\mu}}(\gamma) = \hat{\boldsymbol{\mu}}_k + \gamma(\bar{\mathbf{y}}_{k+1} - \hat{\boldsymbol{\mu}}_k)$.

Note that, by definition, the value of $|\hat{c}_j|$ is the same for all $j \in \mathcal{A}_0$. Hence, the functions (2) for those variables are identical, namely $(1 - \gamma)\hat{C}$, and for all $j \in \mathcal{A}_1$ the corresponding functions $|\hat{c}_j(\gamma)|$ will intersect $(1-\gamma)\hat{C}$ at $\gamma = 1$. This explains why in line 6 we only have to check for the first intersection between $(1 - \gamma)\hat{C}$ and $|\hat{c}_j(\gamma)|$ for $j \in \mathcal{A}^c$.

It also follows from (2) that, for all $j \in \mathcal{A}_0$, we have

$$\mathbf{x}'_j(\bar{\mathbf{y}}_{k+1} - \hat{\boldsymbol{\mu}}_k) = \text{sign}(\hat{c}_j)\hat{C}.$$

Thus, for those variables that are in $\mathcal{A}_0$ we move in line 6 of the modified algorithm in a direction that has a similar geometric interpretation as the direction along which we move in the original LARS algorithm. Namely that for each $j \in \mathcal{A}_0$ the angle between the direction in which we move and $\text{sign}(\hat{c}_j)\mathbf{x}_j$ is the same and this angle is less than $90°$.

Figure 5 shows the result of the modified LARS analysis for the same data used above. Putting variables that enter the model simultaneously into brackets, the order in which the variables enter the model is $(X_{2:5}, X_2, X_5)$, $(X_{3:4}, X_3, X_4)$, $X_{2:5}$, $X_{2:3}$, $(X_{1:1}, X_1), \ldots$. That is, the modified LARS algorithm selects in this case in five steps a model with 10 terms, 6 of which are the terms that are indeed in model (1).

Using the same 1000 samples of size $n = 500$ as above and performing a modified LARS analysis on them using a design matrix with all main and interaction terms also shows markedly improved results. Again, for each replication the step at which a variable was selected into the model was recorded and Figure 6 shows for each variable histograms for these data. To avoid cluttering, the histograms in Figure 6 were truncated to $[1, 20]$; the complete histograms are shown on the right in Figure 7.

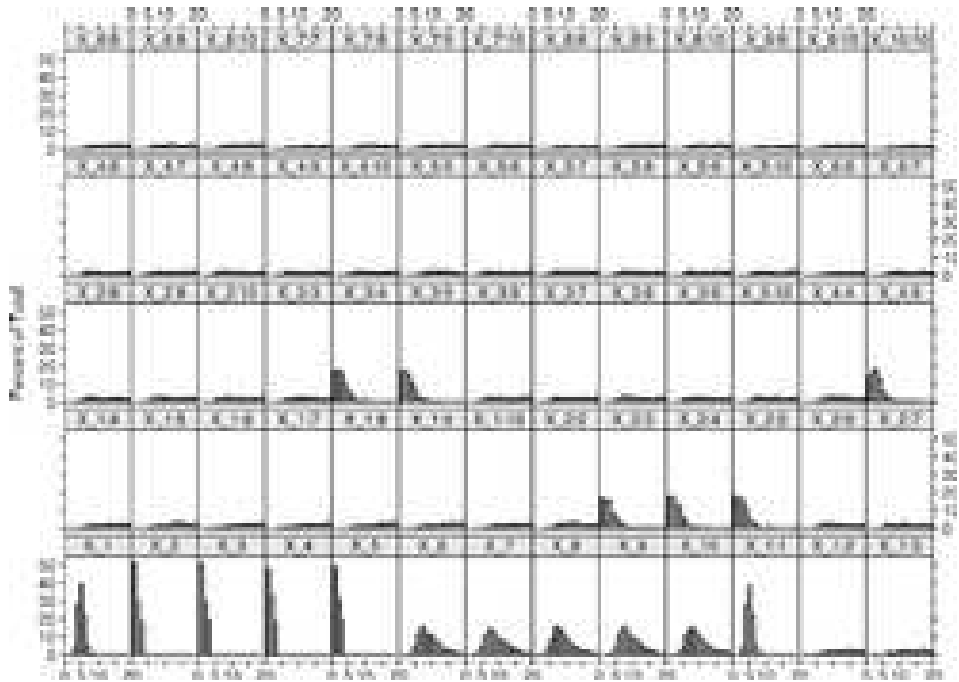

Fig. 6. *Percentage histogram of step at which each variable is selected based on* 1000 *replications: results shown for variables selected in the first* 20 *steps of a modified LARS analysis using main and interaction terms.*

From Figure 6 it can be seen that now the variables $X_2$, $X_3$, $X_4$ and $X_5$ are all selected within the first three iterations of the modified LARS algorithm. Also $X_{1:1}$ and $X_1$ are picked up consistently and early. Compared with Figure 4 there are marked differences in the distribution of when the variable is selected for the interaction terms $X_{i:j}$, $i, j \in \{2, 3, 4, 5\}$, $i \leq j$, and the main terms $X_i$, $i = 6, \ldots, 10$. The latter can be explained by the fact that the algorithm now enforces the marginality principle. Thus, it seems that this modification does improve the performance of the LARS algorithm for model (1). Hopefully it would do so also for other models.

In conclusion, I offer two further remarks and a question. First, note that the modified LARS algorithm may also be used to incorporate factor variables with more than two levels. In such a situation, I would suggest that indicator variables for *all* levels are included in the initial design matrix; but this would be done mainly to easily calculate all the correlations. The dependencies $d_{ij}$ would be set up such that if one indicator variable is selected, then all enter the model. However, to avoid redundancies one would only put all but one of these columns into the matrix $X_{\mathcal{A}}$. This would also avoid that $X_{\mathcal{A}}$ would eventually become singular if more than one explanatory variable is a factor variable.

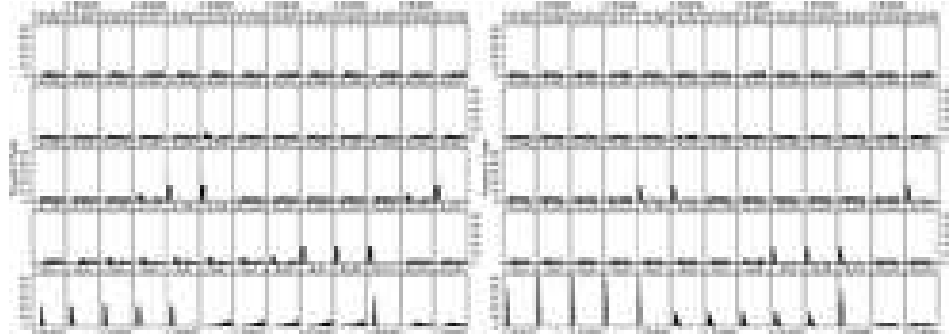

Fig. 7. *Percentage histogram of step at which each variable is selected based on* 1000 *replications:* (left) *LARS analysis;* (right) *modified LARS analysis.*

Second, given the insight between the LARS algorithm and the Lasso algorithm described by EHJT, namely the sign constraint (3.1), it now seems also possible to modify the Lasso algorithm to incorporate the marginality principle by incorporating the sign constraint into Algorithm 2. However, whenever a variable would be dropped from the set $\mathcal{A}_0$ due to violating the sign constraint there might also be variables dropped from $\mathcal{A}_1$. For the latter variables these might introduce discontinuities in the traces of the corresponding parameter estimates, a feature that does not seem to be desirable. Perhaps a better modification of the Lasso algorithm that incorporates the marginality principle can still be found?

Finally, the behavior of the LARS algorithm for model (1) when all main terms and interaction terms are used surprised me a bit. This behavior seems to raise a fundamental question, namely, although we try to build a linear model and, as we teach our students, correlation "measures the direction and strength of the linear relationship between two quantitative variables" [Moore and McCabe (1999)], one has to wonder whether selecting variables using correlation as a criterion is a sound principle? Or should we modify the algorithms to use another criterion?

## REFERENCES

Breiman, L. (1995). Better subset regression using the nonnegative garrote. *Technometrics* **37** 373–384. MR1365720

McCullagh, P. and Nelder, J. A. (1989). *Generalized Linear Models*, 2nd ed. Chapman and Hall, London. MR727836

Moore, D. S. and McCabe, G. P. (1999). *Introduction to the Practice of Statistics*, 3rd ed. Freeman, New York.

Nelder, J. A. (1977). A reformulation of linear models (with discussion). *J. Roy. Statist. Soc. Ser. A* **140** 48–76. MR458743

Nelder, J. A. (1994). The statistics of linear models: Back to basics. *Statist. Comput.* **4** 221–234.

School of Mathematics and Statistics
University of Western Australia
35 Stirling Highway
Crawley WA 6009
Australia
e-mail: berwin@maths.uwa.edu.au